\documentclass{amsart}
\newcommand{\note}{\noindent {\bf Notation. }}
\newcommand{\remark}{\noindent {\bf Remark. }}

\newcommand{\ws}{\hspace{4pt}}
\newtheorem{theorem}{Theorem}
\newtheorem{statement}{Statement}
\newtheorem{proposition}{Proposition}
\newtheorem{cor}{Corollary}
\newtheorem{lemma}{Lemma}

\usepackage[]{amssymb,amsopn}
\usepackage[]{amsmath}
\usepackage[]{eucal}
\usepackage[]{latexsym}

\begin{document}

\title[]{Multiplication operator and exceptional Jacobi polynomials}
\author{\'A. P. Horv\'ath }

\thanks{Supported by the NKFIH-OTKA Grants K128922 and K132097.}

\begin{abstract}
Below the normalized weighted reciprocal of the Christoffel function with respect to exceptional Jacobi polynomials is investigated. It is proved that it tends to the equilibrium measure of the interval of orthogonality in weak-star sense. The main tool of this study is the multiplication operator and examination of behavior of zeros of the corresponding average characteristic polynomial. Finally, as an application of multiplication operator, location of zeros of certain self-inversive polynomials are examined.
\end{abstract}
\maketitle

\section{Introduction}

First we sketch some problems which are strongly related to each other. These problems have an extended literature, below only some examples are cited. \\
Let $\mu$ be a measure on a real interval $I$ with finite moments (or $w>0$ a weight function, where $d\mu(x)=w(x)dx$) . Let $\{q_k\}_{k=0}^\infty$ be the standard orthonormal polynomials on $I$ with respect to $\mu$ (or $w$). Let  $\xi_i:=\xi_{i,k}$, $i=1, \dots, k$ the zeros of $q_k$. In several cases (also in more general circumstances) it is proved that the normalized counting measure based on the zeros of standard orthogonal polynomials tends to the equilibrium measure of the interval of orthogonality is weak-star sense, that is
\begin{equation}\label{lim1}\lim_{n \to \infty}\frac{1}{n}\sum_{i=1}^n\delta_{\xi_i}=\mu_{e,I},\end{equation}
see e.g. \cite{ka} and the references therein.\\
On the other hand the measure, defined as the weighted reciprocal of the Christoffel function, tends to the equilibrium measure of $I$ again:
\begin{equation}\label{lim2}\lim_{n \to \infty}d\mu_n(x)=\lim_{n\to \infty}\frac{1}{n}K_n(x,x)d\mu(x)=\mu_{e,I},\end{equation}
see e.g. \cite{mnt}, \cite{ha1}.\\
The three-term recurrence relation fulfilled by standard orthogonal polynomials,
implies that the multiplication operator, $M: f(x)\to xf(x)$, acting on the weighted $L^2_w$ space has a close relation to the previous two measures. As $M$ can be represented by a tridiagonal (Jacobi) matrix, denoting by $\pi_n$ the projection operator to the $n$-dimensional subspace, the eigenvalues of $\pi_nM\pi_n$ are just the zeros of $q_n$ and limits like \eqref{lim1} and \eqref{lim2} can be proved, see e.g.  \cite{st}, \cite{s}. \\
On the other hand, considering the probability density\\
$\varrho_{n}(x_1, \dots ,x_n)=c({n})\det|K_n(x_i,x_j)|_{i,j=1}^n\prod_{i=1}^n W(x_i),$
where $c(n)$ is a normalization factor, the average characteristic polynomial is defined as\\
$\chi_n(z):=\mathbb{E}\left(\prod_{i=1}^n(z-x_i)\right),$ where the expectation $\mathbb{E}$ refers to $\varrho_n$.
According to \cite[(2.2.10)]{sz}, the zeros of the average characteristic polynomial are just the zeros of $q_n$ again. Thus \eqref{lim1} and \eqref{lim2} can be proved by this technique as well, see e.g. \cite{ha1}, \cite{ha2}.

Below we deal with exceptional Jacobi polynomials. Exceptional orthogonal polynomials are complete systems of polynomials with respect to a positive measure. They are different from the standard or from the classical orthogonal polynomials, since exceptional families have finite codimension in the space of polynomials. Similarly to the classical ones exceptional polynomials are eigenfunctions of Sturm-Liouville-type differential operators but unlike the classical cases, the coefficients of these operators are rational functions. Exceptional orthogonal polynomials also possess a Bochner-type characterization as each family can be derived from one of the classical families applying finitely many Darboux transformations, see \cite{ggm}.\\
Exceptional orthogonal polynomials were introduced recently by G\'omez-Ullate, Kamran and Milson, cf. e.g. \cite{gukm0}, \cite{gukm2} and the references therein. These families of polynomials play a fundamental role for instance in the construction of bound-state solutions to exactly solvable potentials in quantum mechanics. In the last few years have seen a great deal of activity in this area both by mathematicians and physicists, cf. e.g. \cite{d}, \cite{o4}, \cite{hos},  \cite{gugm}, etc. The location of zeros of exceptional orthogonal polynomials is also examined, cf. e.g. \cite{dy}, \cite{gumm}, \cite{km}, \cite{h1}, \cite{h2}, \cite{bo}.

The results on zero-distribution summarized above can be derived considering the three-term recurrence relation fulfilled by standard orthogonal polynomials. Since exceptional orthogonal polynomials fulfil $2L+1$ recurrence formulae with $L\ge 2$, the situation here is different. In \cite{h} by combinatorial methods it is proved that in one codimensional Jacobi case $\mu_n$ tends to $\mu_{e,I}$ in weak-star sense. The relation between the transformed, normalized zero-counting measure based on the zeros of the modified average characteristic polynomial and $\mu_n$ is also studied there.

The aim of this investigation is to extend \eqref{lim2} to exceptional Jacobi polynomials of any codimension. To examine the relation between the (modified) average characteristic polynomial and the exceptional orthogonal polynomials the main tool is the multiplication operator. The location of zeros of exceptional Jacobi polynomials is investigated by outer ratio asymptotics. Finally the multiplication operator method is applied to determine the zeros on the unite circle
of certain self-inversive polynomials.

\section{Notation and the main result}

\subsection{General construction of exceptional orthogonal polynomials.}
There are two different approaches of derivation exceptional orthogonal polynomials from classical ones. Both methods are based on the quasi-rational eigenfunctions of the classical (Hermite, Laguerre, Jacobi) differential operators. The method, developed first for exceptional Hermite polynomials (see \cite{gugm}), defines exceptional polynomials by the Wronskian of certain partition of quasi-rational eigenfunctions mentioned above. The other method defines exceptional polynomials recursively, that is by application of finitely many Crum-Darboux transformations to the classical differential operators. A Bochner-type characterization theorem (see \cite{ggm}) ensures that all exceptional classes can be derived by this method.

Subsequently we follow the construction by Crum-Darboux transformation. According to \cite[Propositions 3.5 and 3.6]{ggm} it is as follows.\\
Classical orthogonal polynomials $\left\{P_n^{[0]}\right\}_{n=0}^\infty$ are eigenfunctions of the second order linear differential operator {with polynomial coefficients}
$$T[y]=py''+qy'+ry,$$
and its eigenvalues are denoted by $\lambda_n$.
$T$ can be decomposed as
\begin{equation}\label{si}T=BA+\tilde{\lambda}, \ws \mbox{with} \ws A[y]=b(y'-wy), \ws B[y]=\hat{b}(y'-\hat{w}y),\end{equation}
where  $w$ is the logarithmic derivative of a quasi-rational eigenfunction of $T$ with eigenvalue $\tilde{\lambda}$, that is $w$ is rational and fulfils
the Riccati equation
\begin{equation}\label{w}p(w'+w^2)+qw+r=\tilde{\lambda}.\end{equation}
$b$ is a suitable rational function.  The coefficients of $B$ can be expressed as
\begin{equation}\label{s0}\hat{b}=\frac{p}{b}, \ws \ws \hat{w}=-w-\frac{q}{p}+\frac{b'}{b}.\end{equation}
Then the exceptional polynomials are the eigenfunctions of $\hat{T}$, that is the partner operator of $T$, which is
\begin{equation}\label{ka}\hat{T}[y]=(AB+\tilde{\lambda})[y]=py''+\hat{q}y'+\hat{r}y,\end{equation}
where
\begin{equation}\label{hat}\hat{q}=q+p'-2\frac{b'}{b}p, \ws \hat{r}=r+q'+wp'-\frac{b'}{b}(q+p')+\left(2\left(\frac{b'}{b}\right)^2-\frac{b''}{b}+2w'\right)p,\end{equation}

\eqref{si} and \eqref{ka} imply that
\begin{equation}\label{sk}\hat{T}AP_n^{[0]}=\lambda_nAP_n^{[0]},\end{equation}
so exceptional polynomials can be obtained from the classical ones by application of (finitely many) appropriate first order differential operator(s) to the classical polynomials. This observation motivates the notation below
\begin{equation}\label{A} AP_n^{[0]}=b\left(P_n^{[0]}\right)'-bwP_n^{[0]}=:P_n^{[1]},\end{equation}
(and recursively $A_sP_n^{[s-1]}=:P_n^{[s]}$ in case of $s$ Darboux transformations.)
The degree of $P_n^{[1]}$ is usually greater than $n$. $\left\{P_n^{[1]}\right\}_{n=0}^\infty$ is an orthogonal system on $I$ with respect to the weight
\begin{equation}\label{ew}W:=\frac{pw_0}{b^2},\end{equation}
where {$w_0$} is one of the classical weights.

\medskip

\remark
- Since at the endpoints of $I$ (if there is any) $p$ may possesses zeros, $b$ can be zero here as well, but $b$ does not have zeros inside $I$,
otherwise the moments of $W$ would not be finite.\\
- As each $P_n^{[1]}$ is a polynomial ($n=0, 1, \dots$), applying the operator $A$ to $P_0^{[0]}$ it  can be seen that $bw$ must be a polynomial and to
$P_1^{[0]}$ shows that $b$ itself is also a polynomial.\\
- Let us recall that $r=0$ in the classical differential operators. Considering \eqref{w} and comparing degrees, $w$ itself can not be a polynomial. By the same reasons $pw$ can not be a polynomial unless it is of degree one. $w$ is a rational function and it has no poles in $(-1,1)$.\\
- Expressing \eqref{si} as
\begin{equation}\label{si1} bBP_n^{[1]}=p\left(P_n^{[1]}\right)'+\left(pw+q-p\frac{b'}{b}\right)P_n^{[1]}=(\lambda_n-\tilde{\lambda})bP_0^{[0]},\end{equation}
it can be easily seen that if $P_n^{[1]}$ had got a double zero at $x_0\in \mathrm{int}I$, then $P_0^{[0]}(x_0)=0$, and by \eqref{A} $\left(P_n^{[0]}\right)'(x_0)=0$,
which is impossible, that is that zeros of $P_n^{[1]}$ which are in the interior of the interval of orthogonality are simple.\\
- Let $a$ be (one of) the finite endpoint(s) of the interval of orthogonality. Again by \eqref{si1} one can derive that if $b(a)\neq 0$,
then $P_n^{[1]}(a)\neq 0$ as well.\\
- If there was an $n\in\mathbb{N}$ such that $P_n^{[1]}(a)= 0$, then $b(a)= 0$ and by \eqref{A} $(bw)(a)=0$ and so $P_n^{[1]}(a)= 0$ for all $n$. That is,
taking $\tilde{P}_n^{[1]} :=b_1\left(P_n^{[0]}\right)'-b_1wP_0^{[0]}$, where $b_1(x)=\frac{b(x)}{x-a}$, we arrive to the exceptional system orthogonal with respect to
$W=\frac{pw_0}{b_1^2}$. Thus we can assume that $\forall$ $n$ $P_n^{[1]}(a)\neq 0$, and if $b(a)=0$, then $(bw)(a)\neq 0$. \\

\subsection{Exceptional Jacobi polynomials.}
Let the $n^{th}$ Jacobi polynomial defined by Rodrigues' formula:
$$(1-x)^{\alpha}(1+x)^{\beta}P_n^{\alpha,\beta}(x)=\frac{(-1)^n}{2^nn!}\left((1-x)^{\alpha+n}(1+x)^{\beta+n}\right)^{(n)},$$
where $\alpha, \beta >-1$.
$$p_k:=p_k^{\alpha,\beta}=\frac{P_k^{\alpha,\beta}}{\varrho_k^{\alpha,\beta}},$$
\begin{equation}\label{ro}\varrho_k^2:=\left(\varrho_k^{\alpha,\beta}\right)^2= \frac{2^{\alpha+\beta+1}\Gamma(k+\alpha+1)\Gamma(k+\beta+1)}{(2k+\alpha+\beta+1)\Gamma(k+1)\Gamma(k+\alpha+\beta+1)}.\end{equation}

the orthonormal Jacobi polynomials which fulfil the following differential equation (cf. \cite[(4.2.1)]{sz})
\begin{equation}\label{jd}(1-x^2)y''+(\beta-\alpha-(\alpha+\beta+2)x)y'+n(n+\alpha+\beta+1)y=0.\end{equation}
With
\begin{equation}\label{jac}P_{n}^{[0]}=p_n^{(\alpha,\beta)}=p_n\end{equation}
and
\begin{equation}\label{ejac}P_{n}^{[1]}=AP_{n}^{[0]}=b(P_{n}^{[0]})'-bwP_{n}^{[0]}.\end{equation}
The next examples of $X_m$ Jacobi polynomials, that is exceptional Jacobi polynomials given by one Darboux transformation and of codimension $m$ can be found in \cite{gumm}. These are as follows.\\
\begin{equation}w_0=w^{(\alpha,\beta)}=(1-x)^\alpha(1+x)^\beta,\end{equation}
where $\alpha$ $\beta$ are defined appropriately, see \cite[Proposition 5.1]{gumm}.
$$T[y]=(1-x^2)y''+(\beta-\alpha-(\alpha+\beta+2)x)y'=BA-(m-\alpha)(m+\beta+1),$$
where $A$ and $B$ are defined in \eqref{si}, and
$$b(x)=(1-x)P_m^{(-\alpha,\beta)}(x), \ws \ws \ws w(x)=(\alpha-m)\frac{P_m^{(-\alpha-1,\beta-1)}(x)}{(1-x)P_m^{(-\alpha,\beta)}(x)}.$$
So the defined exceptional Jacobi polynomials, $P_n^{[1]}:=AP_n^{[0]}$, are orthogonal with respect to
$$W(x)=\frac{(1-x^2)w^{(\alpha,\beta)}(x)}{(1-x)^2\left(P_m^{(-\alpha,\beta)}(x)\right)^2}=\frac{w^{(\alpha-1,\beta+1)}(x)}{\left(P_m^{(-\alpha,\beta)}(x)\right)^2},$$
and the space spanned by these classes are $m$-codimensional in the space of polynomials.

We  restrict our investigations to one-step Darboux transformation case. As it is mentioned above, $b$ can be zero at $\pm 1$. Indeed, by \cite[Table 1]{bo} the quasi-rational eigenfunctions of 
$$T[y]=(1-x^2)y''+(\beta-\alpha-(\alpha+\beta+2)x)y'$$
cf. \eqref{jd}, are $p_n^{(\alpha,\beta)}(x)$, $(1+x)^{-\beta}p_n^{(\alpha,-\beta)}(x)$, $(1-x)^{-\alpha}p_n^{(-\alpha,\beta)}(x)$ and $(1+x)^{-\beta}(1-x)^{-\alpha}p_n^{(-\alpha,-\beta)}(x)$. It gives three kinds of $b(x)$. According to the three types of quasi-rational eigenfunctions let
$$b_0(x)=\left\{\begin{array}{ll}(1-x)p_n^{(-\alpha,\beta)}(x) \ws \ws \mbox{or}\\
(1+x)p_n^{(\alpha,-\beta)}(x) \ws \ws \mbox{or}\\
(1-x^2)p_n^{(-\alpha,-\beta)}(x).\end{array}\right.$$
Thus $b(x)$ can be taken as
\begin{equation}\label{b0}b(x)=s(x)b_0(x),\end{equation}
where $s$ is a polynomial such that $s(x)\neq 0$ on $[-1,1]$.\\
Of course, $\{P_n^{[1]}\}$ is an (exceptional) closed orthogonal system with $s\equiv 1$ in $L^2_{W}$ if and only if $\{sP_n^{[1]}\}$ is orthogonal and closed in $L^2_{\frac{W}{s^2}}$. If $s$ had a zero in $x=\pm 1$, then closedness could fail. So the examples beyond $b_0$ do not give really new classes and can be handled as simple consequences of the $b_0$-type families. Because it does not cause any difficulty, below we use $b$ instead of using $b_0$ first and than taking extension.

Subsequently we denote by $\tilde{b}$ that part of $b$ which appears in the weight function,
\begin{equation}\label{pperb}\frac{p}{b}=\frac{\tilde{p}}{\tilde{b}},\end{equation}
which is bounded on $[-1,1]$.
Introducing the notation
\begin{equation}\label{btilde}b(x)=(1-x)^{\frac{1-\varepsilon_1}{2}}(1+x)^{\frac{1-\varepsilon_2}{2}}\tilde{b}(x),\end{equation}
where $\varepsilon_i=\pm 1$, $i=1,2$;
the exceptional Jacobi polynomial system, $\left\{P_{n}^{[1]}\right\}_{n=0}^{\infty}$, is orthogonal on $(-1,1)$ with respect to
\begin{equation}\label{W} W=\frac{w^{(\alpha+\varepsilon_1,\beta+\varepsilon_2)}}{\tilde{b}^2}.\end{equation}
The orthonormal sytem  is denoted by $\left\{\hat{P}_{n}\right\}_{n=0}^{\infty}.$
The codimension is given by the degree of $\tilde{b}$ cf. \cite{ggm}.

\medskip

\subsection{Formulation of the main theorem}

Before stating the main theorem, we collect the known results about zeros of exceptional Jacobi polynomials. For exceptional Jacopi polynomials defined as the Wronkian of quasi-rational eigenfunctions of the Jacobi operator it is proved that the regular zeros, i.e. zeros in $(-1,1)$ are simple and the normalized zero-counting measure based on regular zeros tends to the equilibrium measure of $[-1,1]$ in weak-star sense, see \cite[Theorem 6.5]{bo}. It is also proved that the exceptional zeros, that is the zeros out of $(-1,1)$ tend to the zeros of the polynomial in the denominator of the weight function, and the speed of convergence is estimated as well, see \cite[Theorem 6.5]{bo}. A Mehler-Heine type asymptotics is also given, see \cite[Theorem 6.3]{bo}. \\
Below these results are derived by Crum-Darboux transformation approach in 1-step case, see Proposition 5, Proposition 3 and Proposition 4, respectively.

The next theorem extends \eqref{lim2} to exceptional Jacobi polynomials of any codimension, that is it generalizes \cite[Theorem 4.1]{h}.

\begin{theorem}\label{T1} Let $\left\{\hat{P}_{n}\right\}_{n=0}^{\infty}$ be the orthonormal system of exceptional Jacobi polynomials generated by one Darboux transformation from the original Jacobi polynomials and of arbitrary codimension. If  $\alpha+\varepsilon_1, \beta+\varepsilon_2 \ge -\frac{1}{2}$, then
\begin{equation}\label{mte}\mu_n \to \mu_e \end{equation}
in weak-star sense, where
$$ d\mu_n(x)=\frac{1}{n}\sum_{k=0}^{n-1}\hat{P}_{k}^2(x)W(x)dx$$
and $\mu_e$ is the equilibrium measure of $[-1,1]$.

\end{theorem}

\medskip

\section{Multiplication operator and proof of Theorem \ref{T1}}

To prove Theorem \ref{T1} we investigate the location of zeros of the modified characteristic polynomial introduced in \cite{h}. In our case it is as follows.\\
Let us denote by $\{\hat{P}_n\}_{n=0}^\infty $ the orthonormal system of exceptional Jacobi polynomials, that is
\begin{equation}\label{hs} \hat{P}_n:=\frac{P_n^{[1]}}{\sigma_{n}}, \ws \ws \mbox{where} \ws \ws \sigma_{n}:=\|P_n^{[1]}\|_{W,2}.\end{equation}
The Christoffel-Darboux kernel is
\begin{equation}K_n(x,y):=\sum_{k=0}^{n-1}\hat{P}_k(x)\hat{P}_k(y).\end{equation}
The joint probability density generated by an $n$-point ensemble has the form
\begin{equation}\label{rho}\rho_{N}(x_1, \dots ,x_N)=c(N)\det|K_N(x_i,x_j)|_{i,j=1}^N\prod_{i=1}^N W(x_i).\end{equation}
Thus the correlation functions are
$$\rho_{N,n}(x_1, \dots ,x_n)=c(n,N)\det|K_N(x_i,x_j)|_{i,j=1}^n\prod_{i=1}^nW(x_i),$$
where $c(N)$, $c(n,N)$ are normalization factors. The expectation $\mathbb{E}$ refers to \eqref{rho}.
The eigenvalue measure becomes a point process with determinantal correlation kernel, cf. \cite[Example 2.12]{j}, i.e. it is a determinaltal point process.

Now we define the modified characteristic polynomial. Let $b$ and $\tilde{b}$ be given by \eqref{ejac} and \eqref{pperb}, respectively. Recalling that $b$ is a polynomial, let
\begin{equation}\label{Q} Q(x):= \int^x \tilde{b}.\end{equation}
Since the definition of $Q$ let the constant term be chosen, we choose it to be zero, say.\\
Considering $Q$, the modified average characteristic polynomial is defined as
\begin{equation}\chi_n(z):=\chi_n^{Q}(z)=\mathbb{E}\left(\prod_{i=1}^n(z-Q(x_i))\right),\end{equation}
cf. \cite{h}. Denote by $z_i$ the zeros of $\chi_n(z)$. Define the normalized zero-counting measure by
$\nu_{{n}}=\frac{1}{n}\sum_{i=1}^n\delta_{z_i},$ and the modified empirical distribution,
$\hat{\mu}_n^{Q}=\frac{1}{n}\sum_{i=1}^n\delta_{Q(x_i)}.$
In 1-step exceptional Jacobi case it is proved that for all $l\ge 0$
\begin{equation}\label{mcp}\lim_{n\to \infty}\left|\mathbb{E}\left(\int x^ld\hat{\mu}_n^{Q}(x)\right)- \int x^ld\nu_n(x)\right|=0,\end{equation}
see \cite[Theorem 5.1]{h}. That is the examination of $\nu_{{n}}$ leads to the description of the behavior of the normalized "Christoffel function measure", i.e.
\begin{equation}d\mu_n(x)=\frac{1}{n}K_n(x,x)W(x)dx.\end{equation}

Indeed, define the not necessarily unique  preimages of $z_i$  by $y_i$, i.e. $z_i=Q(y_i)$, $i=1, \dots , n$. The normalized counting measure based on $y_i$ is
\begin{equation}\label{nuhu}\tilde{\nu}_n=\frac{1}{n}\sum_{i=1}^n\delta_{y_i}\end{equation}
that is $\int Q^l(y)d\tilde{\nu}_n=\int z^ld\nu_n(z)$. Note, that since $\tilde{b}>0$ on $[-1,1]$, $Q$ is increasing here and so if $z_i \in Q([-1,1])$, $y_i$ can be chosen from $[-1,1]$ uniquely. \\
According to \cite[Proposition 2.2]{j} \eqref{mcp} reads
\begin{equation}\label{mz}\lim_{n\to \infty}\left|\int_{-1}^1Q^ld\mu_n- \int Q^ld\tilde{\nu}_n(x)\right|=0,\end{equation}
for all $l\ge 0$, see \cite[Corollary 5.2]{h}.

The next theorem describes zero distribution of the modified average characteristic polynomial. $X_m$ stands for exceptional orthogonal polynomials generated by one Darboux transformation and of codimension $m$.

\medskip

\begin{theorem}\label{nutoe} In $X_m$ exceptional Jacobi case, where $m\in \mathbb{N}$ arbitrary if $\alpha+\varepsilon_1, \beta+\varepsilon_2 \ge -\frac{1}{2}$, then for all $l\in\mathbb{N}$
$$\lim_{n\to\infty}\left(\int Q^l d\tilde{\nu}_n-\int_{-1}^1 Q^l d\mu_{e}\right)=0,$$
where $\mu_e$ is the equilibrium measure of $[-1,1]$.
\end{theorem}

\medskip

For this study our main tool is the multiplication operator generated by $Q$, $M: f \to Qf$, cf. \cite{h}.\\ Exceptional orthogonal polynomials fulfil the next recurrence formula with constant coefficient: $QP_n^{[1]}=\sum_{k=-L}^L\tilde{u}_{n,k}P_{n+k}^{[1]}$ (see  \cite{o} and \cite[(3.4)]{h}).
After normalization the recurrence relation  above is modified as
\begin{equation}\label{R}Q\hat{P}_n=\sum_{k=-L}^Lu_{n,k}\hat{P}_{n+k},\end{equation}
where $u_{n,k}=\frac{\sigma_{n+k}}{\sigma_{n}}\tilde{u}_{n,k}$, for $\sigma_n$ see \eqref{hs}. That is $M$ can be represented by an infinite matrix, $M_{e}$, in the orthonormal basis $\{\hat{P}_n\}_{n=0}^\infty$. This operator, acting on $L^2_W[-1,1]$ and simultaneously on $l_2$ is denoted by $M_e$.
By \eqref{R} the matrix of $M_{e}$ is $2L+1$-diagonal:
\begin{equation}\label{me}M_{e}=\left[\begin{array}{cccccccc}u_{0,0}&u_{0,1}&\dots&\dots&u_{0,L}&0&0&\dots\\u_{1,-1}&u_{1,0}&\dots&\dots&u_{1,L-1}&u_{1,L}&0&\dots\\ \vdots&\vdots&\dots&\ddots&\vdots&\dots&\vdots&\dots\\u_{L,-L}&u_{L,-L+1}&\dots&u_{L,0}&\dots&\dots&u_{L,L}&0\\ 0&u_{L+1,-L}&\dots&\vdots&\dots&\vdots&\dots&u_{L+1,L}\\ \vdots&0&\dots&u_{L+j,-L}& \dots&\vdots&\dots&\dots\end{array}\right].\end{equation}
It can be easily seen that $M_e$ is symmetric since
\begin{equation}\label{ak}u_{k,j}=\int_{-1}^1Q\hat{P}_k\hat{P}_{k+j}W^2=\int_{-1}^1Q\hat{P}_{k+j}\hat{P}_{(k+j)-j}W^2=u_{k+j,-j}.\end{equation}

\medskip

Let us recall that $\{\hat{P}_n\}_{n=0}^\infty$  is an orthonormal system on $[-1,1]$ with respect to $W=\frac{pw_0}{b^2}$, where $w_0$ is a classical Jacobi weight function. Besides this exceptional orthonormal polynomial system, there is the standard orthonormal polynomial system, $\{q_n\}_{n=0}^\infty$, on $[-1,1]$ with respect to $W$. These standard orthonormal polynomials fulfil the three-term recurrence relation
\begin{equation}\label{sr} xq_n=a_{n+1}q_{n+1}+b_nq_n+a_nq_{n-1}\end{equation}
(see e.g. \cite[(3.2.1)]{sz}). Since $W>0$ on $(-1,1)$, by \cite[Theorem 4.5.7]{n} (see also \cite{r}) in formula \eqref{sr} the recurrence coefficients fulfil the asymptotics
\begin{equation}\label{sa} \lim_{n\to\infty}a_n=\frac{1}{2}\ws \ws \mbox{and} \ws \ws \lim_{n\to\infty}b_n=0.\end{equation}
According to \eqref{sr} multiplication operator on $L^2_W[-1,1]$, $A: f \mapsto xf$ can be represented in the (Schauder) basis $\{q_n\}$ as an infinite tridiagonal matrix (denoted by $A$ again)
$$A=\left[\begin{array}{ccccc}b_0&a_1&0&0&\dots\\a_1&b_1&a_2&0&\dots\\0&a_2&b_2&a_3&\dots\\0&0&a_3&b_3&\dots\\ \dots&\dots&\dots&\dots&\dots\end{array}\right],$$
and so $A$ acts on $L^2_W[-1,1]$ and on $l^2$.
Let $Q(x)$ be as above. Then the multiplication operator, $M: f \mapsto Qf$ can be represented as $M=Q(A)$, and in the  basis $\{q_n\}$ it has a matrix $M_q=Q(A)$. Let $\Pi_n^q$ and $\Pi_n^e$ be the projections to $\mathrm{span}\{q_0, \dots, q_{n-1}\}$, and  to  $\mathrm{span}\{\hat{P}_0, \dots, \hat{P}_{n-1}\}$, respectively. Let
$$A_{n\times n}:=\Pi_n^qA\Pi_n^q,  \ws  \ws M_{e,n\times n}:=\Pi_n^eM_e\Pi_n^e, \ws  \ws M_{q,n\times n}:=\Pi_n^qM_q\Pi_n^q.$$

To prove Theorem \ref{nutoe} we compare the trace of the truncated operators in the different bases above.

It is proved (see \cite{h} and the references therein) that the characteristic polynomial of the truncated multiplication matrix $M_{e,n\times n}$, coincides with the modified average characteristic polynomial, that is
\begin{equation}\label{34} \det (zI_n-\pi_n^eM_e\pi_n^e)=\mathbb{E}\left(\prod_{i=1}^n(z-Q(x_i))\right),\end{equation}
where $I_n$ stands for the $n$-dimensional identity operator.
In the standard case the characteristic polynomial of the truncated multiplication matrix $A_{n\times n}$, coincides with the standard average characteristic polynomial, see \cite{ha1}, and the eigenvalues are just the zeros of the standard orthonormal polynomials, $q_n$, see \cite[(2.2.10)]{sz}. Thus we need the next statement.

\medskip

\begin{proposition}\label{pm} With the notation \eqref{W} let us assume that $\alpha+\varepsilon_1, \beta+\varepsilon_2 \ge -\frac{1}{2}$. Then for each $l\in\mathbb{N}$
$$\lim_{n\to\infty}\frac{1}{n}\left(\mathrm{Tr}((Q(A_{n\times n}))^l)-\mathrm{Tr}((M_{e,n\times n})^l)\right)=0.$$
\end{proposition}

\medskip

\noindent Proof of Proposition \ref{pm} is postponed to the next section.

\medskip

\proof (of Theorem \ref{nutoe})
Let $\xi_i=\xi_{i,n}$, $i=1, \dots, n$ the zeros of the standard orthogonal polynomial, $q_n$. Let us recall that $\{\xi_i\}_{i=1}^n$ are the eigenvalues of $A_{n\times n}$.  Since $\mathrm{Tr}((Q(A_{n\times n}))^l)=\frac{1}{n}\sum_{i=1}^nQ^l(\xi_i)$ and by \eqref{34} $\mathrm{Tr}((M_{e,n\times n})^l)=\frac{1}{n}\sum_{i=1}^nQ^l(y_i)$, according to Proposition \ref{pm}
$$\lim_{n\to \infty}\left|\frac{1}{n}\sum_{i=1}^nQ^l(\xi_i)- \int Q^ld\tilde{\nu}_n\right|=0.$$
On the other hand \eqref{sa} implies that the normalized counting measure based on the zeros of the orthogonal polynomials, $q_n$, tends to the equilibrium measure of the interval of orthogonality in weak-star sense (see eg. \cite[Proposition 1.1]{ka} and the references therein). That is
$$\lim_{n\to \infty}\left|\frac{1}{n}\sum_{i=1}^nQ^l(\xi_i)- \int_{-1}^1 Q^ld\mu_e\right|=0.$$
Comparing the two limits the theorem is proved.

\medskip

\proof (of Theorem \ref{T1})
According to \eqref{mz} and Theorem \ref{nutoe}
$$\lim_{n\to \infty}\left(\int_{-1}^1Q^ld\mu_n-\int_{-1}^1 Q^l d\mu_{e}\right)=0$$
for all $l\in\mathbb{N}$.
As it is pointed out in the proof of \cite[Theorem 4.1]{h}, $\mathrm{span}\{Q^l: l\in\mathbb{N}\}$ is dense in $C[-1,1]$, which implies the result.

\medskip

\remark

In the standard case the $n$th average characteristic polynomial coincides with the $n$th orthogonal polynomial.  In exceptional case these two polynomials are different moreover are of different degrees but since the normalized counting measure based on the zeros of exceptional orthogonal polynomials ($\tilde{\mu}_n$) tends to the equilibrium measure in weak-star sense, see \cite[Theorem 6.5]{bo} or Proposition \ref{c2} below, for all $l\in\mathbb{N}$
$$\lim_{n\to\infty} \left(\int Q^ld\tilde{\nu}_n-\int Q^ld\tilde{\mu}_n\right)=0,$$
that is the difference of the normalized zero-counting measures of the exceptional orthogonal and the modified average characteristic polynomials, $\tilde{\nu}_n-\tilde{\mu}_n$, tends to zero in weak-star sense.

\section{ Proof of Proposition \ref{pm}}

First we need some technical lemmas. Let $C=[c_{ij}]_{i,j=0}^\infty$ be an infinite matrix. $c_{ij}= _i C_j$ stands for the elements of $C$. As above its $n${th} principal minor matrix is $C_{n\times n}:=[c_{ij}]_{i,j=0}^{n-1}$. Its $i$th row is $_i C$ and $j$th column is $C_j$.   With this notation we can state the next lemma.

\medskip

\begin{lemma}\label{l1} Let $C=[c_{ij}]_{i,j=0}^\infty$ be a $2k+1$-diagonal, infinite matrix. Let $P=\sum_{l=0}^Mp_lx^l$ be a fixed polynomial. If the entries of $C$ are bounded, that is there is a $K$ such that $|c_{ij}|<K$ for all $0\le i,j <\infty$, then
\begin{equation}\label{tr}\lim_{n\to\infty}\frac{1}{n}\left(\mathrm{Tr}(P(C_{n\times n}))-\mathrm{Tr}((P(C))_{n\times n})\right)=0.\end{equation}
\end{lemma}

\proof
It is obvious that for each $l\in\mathbb{N}$ $C^l$ is $2(kl)+1$ diagonal:\\
Indeed, $c_{ij}=0$ if $|i-j|>k$. Assuming that $(C^l)_{ij}=:C^l_{ij}=0$ if $|i-j|>kl$ we compute $C^{l+1}_{ij}=\langle _iC^l, C_j\rangle=\sum_{|i-p|\le kl \atop |p-j|\le k} C^l_{ip}C_{pj}$. That is  $C^{l+1}_{ij}=0$ $|i-j|>(l+1)k$.

By induction on $l$ we show that the elements of the principal diagonal are coincide in $C_{n\times n}^l$ and $(C^l)_{n\times n}$ except at most the last $k(l-1)$ ones.\\
Let us assume that $(C_{n\times n}^l)_{ij}=((C^l)_{n\times n})_{ij}$ if $i\le n-(l-1)k$ or $j\le n-(l-1)k$. (For $l=1$ it is obviously fulfilled.) Let $i\le n-kl$. By the assumption and by $2kl+1$-diagonality $_i((C_{n\times n})^l)=_i(C^l)$, and similarly if $j\le n-k$, then $(C_{n\times n})_j=C_j$. Thus $((C_{n\times n})^{l+1})_{ij}=C^{l+1}_{ij}$ if $i\le n-kl$ and $j\le n-k$. If $i\le n-kl$, then $C^l_{ip}=0$ if $p>n$  and if $j>n-k$ then $C_j-(C_{n\times n})_j$ starts with $n$ zeros, thus if $i\le n-kl$ and $j> n-k$, $_i((C_{n\times n})^{l+1})_{j}=C^{l+1})_{ij}$ again. For $j\le n-kl$ we have the same by symmetry.\\
Finally considering that the rows and columns contain finitely many non-zero elements, these imply that $\frac{1}{n}\left(\mathrm{Tr}((C_{n\times n})^l)-\mathrm{Tr}((C^l)_{n\times n})\right)\le \frac{1}{n}k(l-1)(2kl+1)^lK^l$, that is
\begin{equation}\label{1sor}\lim _{n\to\infty}\frac{1}{n}\left(\mathrm{Tr}((C_{n\times n})^l)-\mathrm{Tr}(((C^l))_{n\times n})\right)=0.\end{equation}
  Since $\mathrm{Tr}(P(C_{n\times n}))=\sum_{l=0}^Mp_l\mathrm{Tr}((C_{n\times n})^l)$ and $\mathrm{Tr}((P(C))_{n\times n})=\sum_{l=0}^Mp_l\mathrm{Tr}((C^l)_{n\times n})$, \eqref{1sor} implies \eqref{tr}.

\medskip

To ensure the boundedness of the entries of the matrices in question we recall the asymptotic behavior of recurrence coefficients. First we  compute the norm of $P_n^{[1]}$.

\begin{lemma}\label{si} With the notation \eqref{hs}, if $\alpha+\frac{\varepsilon_1}{2}, \beta+\frac{\varepsilon_2}{2}> -\frac{1}{2}$,
\begin{equation}\label{sigma}\sigma_k=\sqrt{k(k+\alpha+\beta+1)+\tilde{\lambda}}.\end{equation}\end{lemma}

\proof  With the notation $p_k^{(\alpha,\beta)}=p_k$, cf. \eqref{jac}
$$\sigma_k^2\hat{P}_k^2W=(bp_k'-bwp_k)^2\frac{pw^{(\alpha,\beta)}}{b^2}$$ $$=(p_k')^2w^{(\alpha+1,\beta+1)}+ (bw)^2p_k^2\frac{pw^{(\alpha,\beta)}}{b^2}-2b^2wp_k'p_k\frac{pw^{(\alpha,\beta)}}{b^2}.$$
By \eqref{w} $pw^2=\tilde{\lambda}-qw-pw'$ and considering \eqref{jd} $(w^{(\alpha+1,\beta+1)})'=qw^{(\alpha,\beta)}$
$$\sigma_k^2\hat{P}_k^2W=(p_k')^2w^{(\alpha+1,\beta+1)}+\tilde{\lambda}p_k^2w^{(\alpha,\beta)}-\left(p_k^2ww^{(\alpha+1,\beta+1)}\right)'.$$
Thus
$$\sigma_k^2=\int_{-1}^1(p_k')^2w^{(\alpha+1,\beta+1)}+\tilde{\lambda}p_k^2w^{(\alpha,\beta)}-\left(p_k^2ww^{(\alpha+1,\beta+1)}\right)'=I_1+I_2+I_3.$$
According to \cite[(4.21.7)]{sz} $p_k'=\frac{k+\alpha+\beta+1}{2}p_{k-1}^{(\alpha+1,\beta+1)}\frac{\varrho_{k-1}^{\alpha+1,\beta+1}}{\varrho_{k}^{\alpha,\beta}}$
$$I_1=\left(\frac{k+\alpha+\beta+1}{2}\frac{\varrho_{k-1}^{\alpha+1,\beta+1}}{\varrho_{k}^{\alpha,\beta}}\right)^2.$$
$I_2=\tilde{\lambda}$.  By the assumption on $\alpha$ and $\beta$ one can see that $I_3=0$. Substituting the values of the corresponding $\varrho_k^{\alpha,\beta}$, \eqref{sigma} is proved.

\medskip

Similarly to \eqref{sa}, the coefficients in \eqref{R} fulfil a symmetric limit relation
\begin{equation}\label{UU}\lim_{n\to\infty}u_{n,j}=:U_{|j|},\end{equation}
where $U_{|j|}$ depends on the polynomial $\tilde{b}$ (cf. \eqref{si}, \eqref{pperb}) as follows. Let
\begin{equation}\label{tb} \tilde{b}(x)=\sum_{k=0}^{L-1} d_kx^k.  \end{equation}
Since Lemma \ref{si} ensures that the asymptotics of $\tilde{u}_{n,j}$ and $u_{n,j}$ coincides, according to \cite[Proposition 1]{h}
\begin{equation}\label{u3}U_{|j|}=\left\{\begin{array}{ll}\sum_{p=\max\{l,1\}}^{\left[\frac{L}{2}\right]}\frac{d_{2p-1}}{2p}\binom{2p}{p-l}\frac{1}{2^{2p}}, \ws\ws \mbox{if}\ws \ws |j|=2l\\
\sum_{p=l}^{\left[\frac{L-1}{2}\right]}\frac{d_{2p}}{2p+1}\binom{2p+1}{p-l}\frac{1}{2^{2p+1}}, \ws\ws \mbox{if}\ws \ws |j|=2l+1. \end{array}\right.\end{equation}

\medskip

Let us recall $\frac{p}{b}$ is bounded on $[-1,1]$ and $W(x)=\frac{(1-x)^{\alpha+\varepsilon_1}(1+x)^{\beta+\varepsilon_2}}{\tilde{b}^2}$ ($\varepsilon_i=\pm 1$), where $b(x)=\tilde{b}(x)(1-x)^{\frac{1-\varepsilon_1}{2}}(1+x)^{\frac{1-\varepsilon_2}{2}}$ and $0<c<\tilde{b}<C$ on $[-1,1]$.\\
Subsequently the next lemma proved by Badkov (see \cite{ba}) is useful. Here we cite the formulation given in \cite{mave}.

\medskip

{\bf Lemma A.}\cite[Lemma 2.E]{mave} {\it Let $\{q_k\}_{k=0}^\infty$ be the standard orthonormal system with respect to $W$. For each $j\ge 0$ integer
\begin{equation}\label{g}\left|q_k^{(j)}(x)\right|\end{equation} $$\le c \left(\frac{k}{\sqrt{1-x^2}+\frac{1}{k}}\right)^j\frac{1}{\left(\sqrt{1-x}+\frac{1}{k}\right)^{\alpha+\varepsilon_1}\left(\sqrt{1+x}+\frac{1}{k}\right)^{\beta+\varepsilon_2}\sqrt{\sqrt{1-x^2}+\frac{1}{k}}}.$$}

\medskip

We need the next estimations on norms of exceptional Jacobi polynomials.

\medskip

\begin{lemma}
If $\alpha+\varepsilon_1, \beta+\varepsilon_2\ge -\frac{1}{2}$, $0\le \delta \le \min\left\{\frac{1}{4},\frac{\alpha+\varepsilon_1}{2}+\frac{1}{4}, \frac{\beta+\varepsilon_2}{2}+\frac{1}{4}\right\}$
\begin{equation}\label{infnorm}\left\|\hat{P}_k(x)\sqrt{W(x)}(1-x^2)^{\frac{1}{4}-\delta}\right\|_\infty \le c k^{\max\{-\varepsilon_1,-\varepsilon_2\}-1+2\delta},\end{equation}
where $c$ is a constant (independent of $k$).
\end{lemma}

\proof
$$\left|\hat{P}_k(x)\sqrt{W(x)}(1-x^2)^{\frac{1}{4}}\right|$$ $$\le \frac{c}{\sigma_k}\left(|b(x)p_k'(x)|+|(bw)(x)p_k(x)|\right)(1-x)^{\frac{\alpha+\varepsilon_1}{2}+\frac{1}{4}}(1+x)^{\frac{\beta+\varepsilon_2}{2}+\frac{1}{4}}$$ $$=K_1+K_2.$$
$$K_1\le \left|p_{k-1}^{(\alpha+1,\beta+1)}\right|(1-x)^{\frac{\alpha}{2}+\frac{3}{4}}(1+x)^{\frac{\beta}{2}+\frac{3}{4}},$$
which is bounded, see e.g. \cite[(8.21.10]{sz}, or \eqref{g}. Since $bw$ is bounded on $[-1,1]$ (it is a polynomial), by \eqref{g}
$$K_2\le \frac{c}{k}|p_k(x)|(1-x)^{\frac{\alpha+\varepsilon_1}{2}+\frac{1}{4}-\delta}(1+x)^{\frac{\beta+\varepsilon_2}{2}+\frac{1}{4}-\delta}$$ $$\le \frac{c}{k}\frac{1}{(\sqrt{1-x}+\frac{1}{k})^\alpha(\sqrt{1+x}+\frac{1}{k})^\beta\sqrt{\sqrt{1-x^2}+\frac{1}{k}}}$$ $$\leq \frac{c\sqrt{k}}{k}\frac{1}{(\sqrt{1-x}+\frac{1}{k})^{-\varepsilon_1-\frac{1}{2}+2\delta}(\sqrt{1+x}+\frac{1}{k})^{-\varepsilon_2-\frac{1}{2}+2\delta}}\le c k^{\max\{-\varepsilon_1,-\varepsilon_2\}-1+2\delta}.$$

\medskip

To prove Proposition \ref{pm} we introduce the operator $O$ which changes the basis in the Hilbert space in question:
\begin{equation}\label{cv} O^{-1}M_qO=M_e, \end{equation}
where
$$O=[o_{ij}]_{i,j=0}^\infty,$$
and
\begin{equation}\label{pq}\hat{P}_j=\sum_{i=0}^{j+m}o_{ij}q_i.\end{equation}

\medskip

\proof (of Proposition \ref{pm})\\
Recalling that $A$ is tridiagonal and $M_e$ is $2L+1$ diagonal and by \eqref{sa} and \eqref{UU} their entries are bounded, by Lemma \ref{l1} it is enough to show that
\begin{equation}\lim_{n\to\infty}\frac{1}{n}\left(\mathrm{Tr}(\pi_n^qQ^l(A)\pi_n^q)-\mathrm{Tr}(\pi_n^e(M_e)^l\pi_n^e)\right)=0.
\end{equation}
Considering \eqref{cv} it is enough to prove that
\begin{equation}\label{36}\lim_{n\to\infty}\frac{1}{n}\left(\mathrm{Tr}(\pi_n^qQ^l(A)\pi_n^q)-\mathrm{Tr}(\pi_n^eO^{-1}Q^l(A)O\pi_n^e)\right)=0.
\end{equation}
To prove \eqref{36} we compute the diagonal elements of $\pi_n^eO^{-1}Q^l(A)O\pi_n^e$. According to \eqref{pq}, for $0\le i\le n-1$
$$_i(O^{-1}Q^l(A)O)_i=\sum_{k=\max\{0,i-L\}}^{\infty}\ _i(O^{-1}Q^l(A))_k \ _kO_i$$
and considering the $2lL+1$-diagonality of $Q^l(A)$
$$_i(O^{-1}Q^l(A))_k =\sum_{j=0}^{i+L}\ _iO^{-1}_j \ _j(Q^l(A))_k=\sum_{\max\{0,k-lL\}\le j \le \min\{i+L, k+lL\}} \ _jO_i \ _j(Q^l(A))_k.$$
Thus
$$\sum_{i=0}^{n-1} \ _i(O^{-1}Q^l(A)O)_i=\sum_{i=0}^{n-1}\sum_{k=0}^{i+L}\sum_{\max\{0,k-lL\}\le j \le \min\{i+L, k+lL\}} o_{ji} o_{ki} \ _j(Q^l(A))_k$$ \begin{equation}\label{hn} =\sum_{i=0}^{n-1}\sum_{k=0}^{i+L}o_{ki}^2\  _k(Q^l(A))_k + \sum_{i=0}^{n-1}\sum_{k=0}^{i+L}\sum_{\max\{0,k-lL\}\le j \le \min\{i+L, k+lL\} \atop j\neq k} o_{ji} o_{ki} \ _j(Q^l(A))_k $$ $$=\tilde{M}_n+H_n.\end{equation}
Reversing the  order of summation
$$\tilde{M}_n=\sum_{k=0}^L\sum_{i=0}^{n-1}o_{ki}^2\  _k(Q^l(A))_k+\sum_{k=L+1}^{n-1}\  _k(Q^l(A))_k\sum_{i=k-L}^{n-1}o_{ki}^2+\sum_{k=n}^{n-1+L}\  _k(Q^l(A))_k\sum_{i=k-L}^{n-1}o_{ki}^2$$ $$=E_{1,n}+M_n+E_{2,n}.$$
Thus the expression after the limit in \eqref{36} becomes
\begin{equation}\label{fosor}\frac{1}{n}\sum_{k=0}^{n-1} \left(\ _k(Q^l(A))_k - \ _k(O^{-1}Q^l(A)O)_k\right)\end{equation}
$$ =\frac{1}{n}\left(\sum_{k=0}^L \ _k(Q^l(A))_k +E_{1,n}+E_{2,n}\right)+\frac{1}{n}\sum_{k=L+1}^{n-1} \ _k(Q^l(A))_k\left(1-\sum_{i=k-L}^{n-1}o_{ki}^2\right) + \frac{1}{n}H_n.$$
According to \eqref{u3} the entries $_i(Q^l(A))_k$ are bounded, and by orthonormality
\begin{equation}\label{oosz}\sum_{i=k-L}^{n-1}o_{ki}^2=1.\end{equation}
Thus
$$\frac{1}{n}\left(\sum_{k=0}^L \ _k(Q^l(A))_k +E_{1,n}+E_{2,n}\right)\le \frac{K}{n},$$
where $K$ depends on $L$ and on the supremum of $|_i(Q^l(A))_k|$, but it is independent of $n$, that is the first term tends to zero, when $n$ tends to infinity.

We can start the estimation of the second term of \eqref{fosor} in the same way. By the boundedness of the entries of $Q^l(A)$, and by \eqref{oosz}
\begin{equation}\label{50}\frac{1}{n}\sum_{k=L+1}^{n-1} \ _k(Q^l(A))_k\left(1-\sum_{i=k-L}^{n-1}o_{ki}^2\right)\le C\frac{1}{n}\sum_{k=0}^{n-1}\sum_{i=n}^\infty o_{ki}^2.\end{equation}
Let $i>k$. Referring to \eqref{pq}
\begin{equation}\label{oki}o_{ki}=\int_{-1}^1\hat{P}_iq_kW=\frac{1}{\lambda_i}\int_{-1}^1\left(p\hat{P}''_i+ \hat{q}\hat{P}'_i+\hat{r}\hat{P}_i\right)q_kW,\end{equation}
where the differential equation of $\hat{P}_i$ was taken into consideration, cf. \eqref{ka} and \eqref{sk}.
By the assumption on $\alpha, \beta$ and \eqref{pperb}, integrating by parts the first term we have
\begin{equation}\label{der}\int_{-1}^1p\hat{P}''_iq_kW=-\int_{-1}^1\hat{P}'_i(pq_kW)'=-\int_{-1}^1\hat{P}'_i(\hat{q}q_k+pq'_k)W,\end{equation}
because $pW'=(\hat{q}-p')W$, cf. \eqref{hat}.\\
So \eqref{oki} and \eqref{der} imply that
\begin{equation}\label{oki1}o_{ki}=\frac{1}{\lambda_i}\int_{-1}^1 -p\hat{P}'_iq'_kW+\hat{r}\hat{P}_iq_kW=I_{ki}+II_{ki}.\end{equation}
Considering \eqref{si1}
\begin{equation}\label{pp}p\hat{P}'_i=\frac{\lambda_i-\tilde{\lambda}}{\sigma_i}bp_i-\left(pw+q-p\frac{b'}{b}\right)\hat{P}_i.\end{equation}
By \eqref{pp}
\begin{equation}\label{1i}I_{ki}=-\frac{\lambda_i-\tilde{\lambda}}{\sigma_i \lambda_i}\int_{-1}^1 bp_iq'_kW+\frac{1}{\lambda_i}\int_{-1}^1\left(pw+q-p\frac{b'}{b}\right)\hat{P}_iq'_kW=I_{ki}^{(1)}+I_{ki}^{(2)}.\end{equation}
By \eqref{pperb} $\left|pw+q-p\frac{b'}{b}\right|$ is bounded on $[-1,1]$. Indeed, only the first term needs some investigation; recalling that $bw$ is a polynomial, notice that $pw=\frac{p}{b}bw$. \\
Thus, according to \eqref{g} and \eqref{infnorm}
$$|I_{ki}^{(2)}|\le \frac{ck}{i^2}\left\|\hat{P}_k(x)\sqrt{W(x)}(1-x^2)^{\frac{1}{4}}\right\|_\infty$$ $$\times\int_{-1}^1\frac{\sqrt{1-x}^{\alpha+\varepsilon_1}\sqrt{1+x}^{\beta+\varepsilon_2}}{(1-x^2)^{\frac{1}{4}}\left(\sqrt{1-x^2}+\frac{1}{k}\right)^\frac{3}{2}\left(\sqrt{1-x}+\frac{1}{k}\right)^{\alpha+\varepsilon_1}\left(\sqrt{1+x}+\frac{1}{k}\right)^{\beta+\varepsilon_2}}dx$$ $$\le c\frac{k}{i^2}J_k.$$
$$J_k\le c \int_{-1}^0\frac{\left(\sqrt{1+x}\right)^{\beta+\varepsilon_2+\frac{3}{2}-2\delta}}{\left(\sqrt{1+x}+\frac{1}{k}\right)^{\beta+\varepsilon_2+\frac{3}{2}-2\delta+2\delta}(1+x)^{1-\delta}}$$ $$+c\int_0^1\frac{\left(\sqrt{1-x}\right)^{\alpha+\varepsilon_1+\frac{3}{2}-2\delta}}{\left(\sqrt{1-x}+\frac{1}{k}\right)^{\alpha+\varepsilon_1+\frac{3}{2}-2\delta+2\delta}(1-x)^{1-\delta}}\le c \frac{k^{2\delta}}{\delta},$$
where the last inequality fulfils if $\alpha+\varepsilon_1+\frac{3}{2}-2\delta\ge 0$ and $\beta+\varepsilon_2+\frac{3}{2}-2\delta\ge 0$.
Let $\delta=ck^{-\frac{1}{4}}$. Then
\begin{equation}\label{ik2}|I_{ki}^{(2)}|\le c\frac{k^{\frac{5}{4}}}{i^2},\end{equation}
provided that $\alpha+\varepsilon_1,\beta+\varepsilon_2> - \frac{3}{2}$.
For sake of simplicity let us denote by $c(i):=-\frac{\lambda_i-\tilde{\lambda}}{\sigma_i \lambda_i}$. After simplification
$$I_{ki}^{(1)}= c(i)\int_{-1}^1\frac{\tilde{p}}{\tilde{b}}q'_kp_iw^{(\alpha,\beta)}.$$
Recaling that $\frac{1}{\tilde{b}}$ is bounded on $[-1,1]$, let $S_{i-k-2}$ be its uniformly best approximating polynomial on $[-1,1]$ of degree $i-k-2$. Then, as its degree is less than $i$,  $S_{i-k-2}\tilde{p}q'_k$ is orthogonal to $p_i$ with respect to $w^{(\alpha,\beta)}$. Thus
$$I_{ki}^{(1)}= c(i)\int_{-1}^1\left(\frac{1}{\tilde{b}}-S_{i-k-2}\right)\tilde{p}q'_kp_iw^{(\alpha,\beta)}.$$
Taking into account \eqref{jd}, \eqref{sigma} and \eqref{si} $c(i) \le c \frac{1}{i}$. Since $\left(\frac{1}{\tilde{b}}\right)'$ is bounded on $[-1,1]$ too, according to the classical Jackson's theorem
$$|I_{ki}^{(1)}|\le  \frac{c}{i(i-k)}\int_{-1}^1\tilde{p}|q'_k||p_i|w^{(\alpha,\beta)}$$ $$\le \frac{c}{i(i-k)}\left\|q'_k(x)\sqrt{W(x)}(1-x^2)^{\frac{3}{4}}\right\|_\infty\int_{-1}^1|p_i(x)|\sqrt{W(x)}\frac{b(x)}{(1-x^2)^{\frac{3}{4}}}dx$$ $$\le c\frac{k}{i(i-k)}\int_{-1}^1|p_i(x)|(1-x)^{\frac{\alpha}{2}-\frac{1}{4}}(1+x)^{\frac{\beta}{2}-\frac{1}{4}},$$
where the norm of $q'_k$ is estimated by \eqref{g}. Finally by \cite[(7.34.1)]{sz} the last integral can be estimated by a constant independently of $i$, that is
\begin{equation}\label{ik1}|I_{ki}^{(1)}|\le  c\frac{k}{i(i-k)}.\end{equation}
 To estimate $II_{ki}$ first we remark that $\hat{r}$ is bounded on $[-1,1]$. Indeed, it is clear that we have to deal with only the endpoints. Recalling that $\frac{p}{b}$ is bounded on $[-1,1]$ and considering \eqref{hat}, $\hat{q}$ is bounded there too. As it is mentioned in the Remark of section 2 $\hat{P}_n(1)\neq 0$, that is due to \eqref{ka} and \eqref{sk} $\hat{r}$ must be bounded in $1$. In $-1$ it is the same. Thus, by Cauchy-Schwarz inequality
\begin{equation}\label{iiki}|II_{ki}|\le\frac{c}{i^2}\int_{-1}^1|\hat{P}_i|\sqrt{W}|q_k|\sqrt{W}\le \frac{c}{i^2}.\end{equation}
Since $\sum_{i=n}^\infty o_{ki}^2\le 1$, according to \eqref{ik2}, \eqref{ik1} and \eqref{iiki}
$$\frac{1}{n}\sum_{k=0}^{n-1}\sum_{i=n}^\infty o_{ki}^2 \le c \frac{1}{n}\sum_{k=n^{\frac{3}{4}}}^{n-n^{\frac{3}{4}}}\sum_{i=n}^\infty o_{ki}^2+ O\left(\frac{1}{n^{\frac{1}{4}}}\right)$$   \begin{equation}\label{fb}\le \frac{c}{n}\sum_{k=n^{\frac{3}{4}}}^{n-n^{\frac{3}{4}}}\sum_{i=n}^\infty\left(\frac{k^{\frac{5}{2}}}{i^4}+ \frac{k^2}{i^2n^{\frac{3}{2}}}\right)+O\left(\frac{1}{n^{\frac{1}{4}}}\right)=O\left(\frac{1}{\sqrt{n}}\right)+ O\left(\frac{1}{n^{\frac{1}{4}}}\right),\end{equation}
so the expression in \eqref{50} tends to zero.

Finally we estimate the error term $\frac{1}{n}H_n$, cf. \eqref{hn} and \eqref{fosor}. Using the properties of $|_j(Q^l(A))_k |$ and $o_{ki}$ again, it can be estimated as
$$\frac{1}{n}|H_n|\le c\frac{1}{n}\sum_{i=0}^{(l-1)L}\sum_{k=0}^{i+L}\sum_{\max\{0,k-lL\}\le j \le i+L \atop j\neq k} |o_{ji}|| o_{ki}|$$ $$+\frac{1}{n}\left|\sum_{i=(l-1)L+1}^{n-1}\sum_{k=0}^{i+L}\sum_{\max\{0,k-lL\}\le j \le  k+lL \atop j\neq k} o_{ji} o_{ki} \ _j(Q^l(A))_k \right|=\Sigma_1(n)+\Sigma_2(n).$$
Certainly, $\lim_{n\to \infty}\Sigma_1(n)=0$. Splitting $\Sigma_2(n)$ to two parts and changing the ordering of summation,
$$\Sigma_2(n)\le\frac{1}{n}\left|\sum_{i=(l-1)L+1}^{n-1}\sum_{k=0}^{lL-1}\sum_{0\le j \le k+lL\atop j\neq k} o_{ji} o_{ki} \ _j(Q^l(A))_k \right| $$ $$+\frac{1}{n}\left|\sum_{i=(l-1)L+1}^{n-1}\sum_{k=lL}^{i+L}\sum_{k-lL\le j \le k+lL \atop j\neq k} o_{ji} o_{ki} \ _j(Q^l(A))_k \right|$$
$$\le\frac{1}{n}\sum_{k=0}^{lL-1}\sum_{0\le j \le k+lL \atop j\neq k}\left| _j(Q^l(A))_k\right|\sum_{i=(l-1)L+1}^{n-1}|o_{ji} o_{ki}|$$ $$+\left|\frac{1}{n}\sum_{k=lL}^{n-1+L}\sum_{k-L\le j \le k+lL \atop j\neq k} \ _j(Q^l(A))_k\sum_{i=(l-1)L+1}^{n-1}o_{ji} o_{ki}\right|=\Sigma_{21}(n)+\Sigma_{22}(n).$$
Again, $\lim_{n\to \infty}\Sigma_{21}(n)=0$. According to \eqref{pq} $o_{ki}=0$ if $i\le k-L$. Thus, by orthogonality $\sum_{i=(l-1)L+1}^{n-1}o_{ji} o_{ki}=\sum_{i=n}^{\infty}o_{ji} o_{ki}$. That is
$$\Sigma_{22}(n)\le c\frac{1}{n}\sum_{k=lL}^{n-1+L}\sum_{i=n}^{\infty}|o_{ji} o_{ki}|\le c\frac{1}{n}\sum_{k=lL}^{n-1+L}\max_{k-L\le j \le k+lL}\sum_{i=n}^{\infty}o_{ji}^2.$$
We can proceed as in \eqref{fb} again, and so  $\lim_{n\to \infty}\Sigma_{22}(n)=0$, which finishes the estimation of $\frac{1}{n}|H_n|$ and the proof of Proposition \ref{pm} as well.

\medskip

\section{Outer ratio asymptotics}

In \cite{gumm} a family of exceptional Jacobi polynomials is given. Among other properties location of their zeros is described there. These results were extended in \cite{bo} by Wronskian method. Below we describe the behavior of zeros and give outer ratio asymptotics and a Heine-Mehler type formula by the general formulation of X$_m$ exceptional Jacobi polynomials. Although some of the results of this section are known, the approach below is different from the one cited above.

\medskip

Recalling the general construction of exceptional orthogonal polynomials with one-step Darboux transformation, cf. \eqref{A} and \eqref{si},
\begin{equation}\label{p}\frac{P_{n-1}^{[1]}}{P_n^{[1]}}=\frac{b(P_{n-1}^{[0]})'-bwP_{n-1}^{[0]}}{b(P_n^{[0]})'-bwP_n^{[0]}}= \frac{P_{n-1}^{[0]}}{P_n^{[0]}}\frac{b\frac{(P_{n-1}^{[0]})'}{P_{n-1}^{[0]}}-bw\frac{P_{n-1}^{[0]}}{P_{n-1}^{[0]}}}{b\frac{(P_n^{[0]})'}{P_n^{[0]}}- bw\frac{P_n^{[0]}}{P_n^{[0]}}}.\end{equation}

Classical orthonormal polynomials satisfy the following relation (cf. e.g. \cite[(4.5.5), (5.1.14), (5.5.10)]{sz}).
\begin{equation}\label{kl}p\left(P_n^{[0]}\right)'=A_nP_{n+1}^{[0]}+B_nP_n^{[0]}+C_nP_{n-1}^{[0]},\end{equation}
where $A_n$, $B_n$ and $C_n$ are {real numbers} and $p$ is the coefficient of the second derivative in the differential equation of the classical orthogonal polynomials. Thus
$$\frac{P_{n-1}^{[1]}}{P_n^{[1]}}=\frac{P_{n-1}^{[0]}}{P_n^{[0]}}\frac{b\left(A_{n-1}\frac{P_{n}^{[0]}}{P_{n-1}^{[0]}}+B_{n-1}+C_{n-1}\frac{P_{n-2}^{[0]}}{P_{n-1}^{[0]}}\right)- bwp}{b\left(A_{n}\frac{P_{n+1}^{[0]}}{P_{n}^{[0]}}+B_{n}+C_{n}\frac{P_{n-1}^{[0]}}{P_{n}^{[0]}}\right)-bwp}.$$

\medskip

Now we return to exceptional Jacobi polynomials.

\medskip

\note Denote by $Z_{\tilde{b}}$ the zeros of the polynomial $\tilde{b}$, cf. \eqref{ejac}, \eqref{pperb}. Similarly to \cite[Proposition 5.6]{gumm} the next asymptotics fulfils.

\medskip

\begin{proposition} For exceptional Jacobi polynomials given by \eqref{ejac}
\begin{equation}\lim_{n\to \infty}\frac{P_{n-1}^{[1]}(z)}{P_n^{[1]}(z)}=z-\sqrt{z^2-1}, \ws \ws z\in\mathbb{C}\setminus [-1,1]\setminus Z_{\tilde{b}}\end{equation}
locally uniformly, where $\sqrt{z^2-1}$ means that branch of the function for which \\$\left|z-\sqrt{z^2-1}\right|<1$ on $\mathbb{C}\setminus [-1,1]$.

\end{proposition}

\proof In Jacobi case \eqref{p} becomes
\begin{equation}\label{pj}\frac{P_{n-1}^{[1]}}{P_n^{[1]}}=\frac{P_{n-1}^{[0]}}{P_n^{[0]}}\frac{n-1}{n}\frac{b\left(\frac{A_{n-1}}{n-1}\frac{P_{n}^{[0]}}{P_{n-1}^{[0]}}+\frac{B_{n-1}}{n-1}+\frac{C_{n-1}}{n-1}\frac{P_{n-2}^{[0]}}{P_{n-1}^{[0]}}\right)-\frac{bwp}{n-1}}{b\left(\frac{A_{n}}{n}\frac{P_{n+1}^{[0]}}{P_{n}^{[0]}}+\frac{B_{n}}{n}+\frac{C_{n}}{n}\frac{P_{n-1}^{[0]}}{P_{n}^{[0]}}\right)-\frac{bwp}{n}}.\end{equation}
Classical Jacobi polynomials fulfil the asymptotics below.
\begin{equation}\label{cj}\lim_{n\to \infty}\frac{P_{n-1}^{[0]}(z)}{P_n^{[0]}(z)}=z-\sqrt{z^2-1}\end{equation}
uniformly on the compact subsets of $\mathbb{C}\setminus [-1,1]$ (cf. e.g. \cite{mnv}).
Taking into consideration that $\lim_{n\to\infty}\frac{\varrho_{n-1}}{\varrho_n}=1$, etc. (for $\varrho_n$ see \eqref{ro}), according to \cite[(4.5.5)]{sz}
\begin{equation}\label{ABC} \lim_{n\to\infty}\frac{A_n}{n}=\frac{1}{2}, \ws \ws \lim_{n\to\infty}B_n = \frac{\alpha-\beta}{2} \ws \ws \lim_{n\to\infty}\frac{C_n}{n}= -\frac{1}{2}.\end{equation}
Since $bwp$ is a polynomial (and usually $pw$ is not) $\frac{bwp}{n}$ tends to zero locally uniformly on $\mathbb{C}$, considering \eqref{cj} and \eqref{ABC} we arrive to the statement.

\medskip

\note
Let us denote by
$$Z:=\{z\in \mathbb{C} :$$ $$ \forall \ws U \ws \mbox{neighborhood of} \ws z \ws \exists \ws N\in \mathbb{N}\ws \mbox{such that} \ws \forall P_n^{[1]} \ws \exists \ws w\in U, \ws P_n^{[1]}(w)=0 \ws \mbox{if} \ws n>N\}.$$

\medskip

\begin{cor}\label{c1} $[-1,1]\subset Z$.
\end{cor}

\proof
For sake of self-containedness we repeat the proof of \cite[Theorem 5]{mnv}.\\ $Z$ is closed. If there was a point $t\in (-1,1)$ which is not in $Z$, there would be a compact neighborhood, $E\subset \mathbb{C}$, of $t$  and a subsequence of polynomials such that $\frac{P_{n_k-1}^{[1]}(z)}{P_{n_k}^{[1]}(z)} \to z-\sqrt{z^2-1}$ on $E$ uniformly, which is impossible.

\medskip

Let us recall that $m=\deg\tilde{b}$ is the codimension of the exceptional system. With this technic in 1-step Darboux-transformation case we get a version of \cite[Theorem 6.6]{bo}.

\medskip

\begin{proposition}\label{exgy}If $n$ is large enough, $P_n^{[1]}$ has $m$ exceptional zeros (with multiplicity), that is $m$ zeros out of the interval of orthogonality. Moreover the exceptional zeros tend to the zeros of $\tilde{b}$, when $n$ tends to infinity.
\end{proposition}

\proof Similarly to the computation above
$$\frac{P_{n}^{[1]}}{nP_n^{[0]}}=\frac{pb(P_n^{[0]})'}{pnP_n^{[0]}}-\frac{bw}{n}=\frac{b}{p}\left(\frac{A_{n}}{n}\frac{P_{n+1}^{[0]}}{P_{n}^{[0]}}+\frac{B_{n}}{n}+\frac{C_{n}}{n}\frac{P_{n-1}^{[0]}}{P_{n}^{[0]}}\right)-\frac{bw}{n},$$
that is by \eqref{cj} and \eqref{ABC}
$$\lim_{n\to\infty}\frac{P_{n}^{[1]}}{nP_n^{[0]}}=-\frac{b(z)}{\sqrt{z^2-1}},$$
where the convergence is locally uniform on $\mathbb{C}\setminus [-1,1]$. Applying Hurwitz's theorem the statement is proved.

\medskip

Let $j_{\alpha}(z)=\Gamma(\alpha+1) \left(\frac{2}{z}\right)^\alpha J_\alpha(z)= \sum_{k=0}^\infty\frac{(-1)^k\Gamma(\alpha+1)}{\Gamma(k+1)\Gamma(k+\alpha+1)}\left(\frac{z}{2}\right)^{2k}$ be the Bessel function.
Classical Jacobi polynomials fulfil the next Mehler-Heine formula (see \cite[Theorem 8.1.1]{sz})
\begin{equation}\label{mh} \lim_{n\to \infty}\frac{P_{n}^{\alpha,\beta}\left(\cos \frac{z}{n}\right)}{n^\alpha}=\frac{1}{\Gamma(\alpha+1)}j_{\alpha}(z),\end{equation}
where the convergence is locally uniform on the complex plane.\\
If $b(1)=0$, that is $\varepsilon_1=-1$ cf. \eqref{btilde}, we introduce the notation $b(x)=(1-x)b_1(x)$. Similarly to \cite[Proposition 5.7]{gumm}, again a 1-step Darboux-transformation case of \cite[Theorem 6.3]{bo} can be derived.

\medskip

\begin{proposition}\label{mehe} If $\alpha\ge -\frac{\varepsilon_1}{2}$
$$\lim_{n\to \infty}\frac{\varrho_n}{n^{\alpha+1+\varepsilon_1}}P_n^{[1]}\left(\cos \frac{z}{n}\right)=cj_{\alpha+\varepsilon_1}(z),$$
where the convergence is locally uniform on the complex plane and $c$ is a constant depending on $\alpha$ and $b$.
\end{proposition}

\medskip

\proof
$\varepsilon_1=1:$\\ By the classical formula
$$(P_{n}^{\alpha,\beta})'=\frac{n+\alpha+\beta+1}{2}P_{n-1}^{\alpha+1,\beta+1},$$
see \cite[(4.21.7)]{sz} and considering \eqref{ejac}
$$\frac{\varrho_n}{n^{\alpha+2}}P_n^{[1]}\left(\cos \frac{z}{n}\right)$$ $$=\frac{n+\alpha+\beta+1}{2n}b\left(\cos \frac{z}{n}\right)
\frac{P_{n-1}^{\alpha+1,\beta+1}\left(\cos \frac{z}{n}\right)}{n^{\alpha+1}}-\frac{(bw)\left(\cos \frac{z}{n}\right)}{n^2}
\frac{P_{n}^{\alpha,\beta}\left(\cos \frac{z}{n}\right)}{n^\alpha}.$$
Since $\|j_{\alpha}\|\le 1$ (see \cite[Ch. 7 7.3 (4)]{be}), applying \eqref{mh} the statement is proved in the first case with $c=\frac{b(1)}{2\Gamma(\alpha+2)}$.\\
$\varepsilon_1=-1:$\\
Taking into consideration the next classical formula (see eg. \cite[(73)]{gumm}):
$$(1-x)(P_{n}^{\alpha,\beta})'=\alpha P_{n}^{\alpha,\beta}-(n+\alpha)P_{n}^{\alpha-1,\beta+1}$$
we have
$$\varrho_n P_n^{[1]}=b_1\left(\alpha
P_{n}^{\alpha,\beta}-(n+\alpha)P_{n}^{\alpha-1,\beta+1}\right)-bwP_{n}^{\alpha,\beta}.$$

Recalling the Remark in Section 2 $(bw)(1)\ne 0$. We show that the polynomial $\alpha b_1-bw$ vanishes at $1$. Indeed, by \eqref{ejac} and \eqref{si}
$$bBP_n^{[1]}(1)=\left(p\left(P_n^{[1]}\right)'+\left(pw+q-p\frac{b'}{b}\right)P_n^{[1]}\right)(1)$$ $$=\left(pw+q-p\frac{b'}{b}\right)(1)(-bwp_n)(1)
=\left(\lambda_n-\tilde{\lambda})b p_n\right)(1)=0.$$
Thus referring to \eqref{jd}
$$0=\left(pw+q-p\frac{b'}{b}\right)(1)=\left(\frac{p}{b}bw+q-p\frac{b'}{b}\right)(1)$$ $$=\frac{1+1}{b_1(1)}(bw)(1)-2\alpha-2-\frac{1+1}{b_1(1)}
(-b_1(1)+(1-1)b_1'(1))=2\left(\frac{(bw)(1)}{b_1(1)}-\alpha\right).$$
That is
$$\varrho_n P_n^{[1]}(x)=(\alpha b_1-bw)(x)P_{n}^{\alpha,\beta}(x)
-b_1(x)(n+\alpha)P_{n}^{\alpha-1,\beta+1}(x)$$ $$=(1-x)s(x)P_{n}^{\alpha,\beta}(x)
-b_1(x)(n+\alpha)P_{n}^{\alpha-1,\beta+1}(x),$$
where $s(x)$ is a polynomial. Applying \cite[(4.5.4)]{sz}
$$\varrho_n P_n^{[1]}=s\frac{2}{2n+\alpha+\beta+1}\left((n+\alpha)P_{n}^{\alpha-1,\beta}-(n+1)P_{n}^{\alpha-1,\beta}\right)
-b_1(n+\alpha)P_{n}^{\alpha-1,\beta+1}.$$
So again by \eqref{mh} and the uniform boundedness of the corresponding Bessel functions
$$\lim_{n\to \infty}\frac{\varrho_n}{n^{\alpha}}P_n^{[1]}\left(\cos \frac{z}{n}\right)=-\frac{b_1(1)}{\Gamma(\alpha)}j_{\alpha-1}(z).$$

\medskip

According to Proposition \ref{exgy} $P_n^{[1]}$ has $m$ exceptional zeros out the interval of orthogonality, and by the Remark of Section 2 $P_n^{[1]}(-1)\neq 0$, $P_n^{[1]}(1)\neq 0$. Because $m$ is the number of gaps in the sequence of degrees, if $n$ is large enough, $P_n^{[1]}$ has to possess $n$ zeros in $(-1,1)$. These are the regular zeros, and as it is noted these zeros are simple. The distribution of regular zeros can be derived from a theorem of Erd\H os and Tur\'an, see \cite{et}.

\medskip

{\bf Theorem B. }\cite{et} {\it  Let $1\ge\zeta_{1,n}>\dots >\zeta_{n,n}\ge-1$, $n\in\mathbb{N}_+$ any system of points, and let $\eta_{i,n}\in[0,\pi]$ be defined by $\zeta_{i,n}=\cos \eta_{i,n}$. Let $\omega_n(\zeta)=\prod_{i=1}^n(\zeta-\zeta_{i,n})$. If  for all $\zeta \in[-1,1]$
$$|\omega_n(\zeta)|\le \frac{A(n)}{2^n}$$
holds, then for every subinterval $[\gamma,\delta] \subset [0,\pi]$ we have
\begin{equation}\label{ert}\left|\sum_{i \atop \gamma\le \eta_{in}\le\delta} 1-\frac{\delta-\gamma}{\pi}n\right|<\frac{8}{\log 3}\sqrt{n\log A(n)}.\end{equation}}

\medskip

A simple application of the previous theorem is the next one.

\medskip

\begin{proposition}\label{c2} Let $x_{1n}, \dots, x_{nn}$ be the regular zeros of the exceptional Jacobi polynomials, $P_n^{[1]}=P_n^{[1],\alpha,\beta}$, $\alpha, \beta \ge -\frac{1}{2}$. Let $x_{in}=\cos\varphi_{in}$.
For every $[\gamma,\delta] \subset [0,\pi]$
$$\left|\frac{1}{n}\sum_{i \atop \gamma\le \varphi_{in}\le\delta} 1-\frac{\delta-\gamma}{\pi}\right|\le c\sqrt{\frac{\log n}{n}},$$
where $c$ is a constant, depends on $\alpha$, $\beta$ $b$, $bw$, but is independent of $n$.
\end{proposition}

\proof
Let $P_n^{[1]}=l_nq_{m,n}s_n$, where $q_{m,n}$ and $s_n$ are monic polynomials of degree $m$ and $n$, respectively, the zeros of $q_{m,n}$ are the exceptional zeros of $P_n^{[1]}$, and $s_n$ possesses the regular ones. $l_n$ is the leading coefficient of $P_n^{[1]}$. Since the zeros of  $q_{m,n}$ tends to the zeros of $\tilde{b}$, if $n$ is large enough, there are constants $k$ and $K$ independent of $n$ such that $0<k<|q_{m,n}|<K$ on $[-1,1]$. Thus for $x\in [-1,1]$
$$|s_n(x)| \le \frac{1}{k|l_n|}|P_n^{[1]}(x)|\le \frac{1}{k|l_n|}(\|bP_n'\|+\|bwP_n\|),$$
where the norm sign refers to the sup-norm on $[-1,1]$, and $P_n=P_n^{(\alpha,\beta)}$. By \cite[(7.32.2)]{sz}, $\|P_n\|=\binom{n+\max\{\alpha,\beta\}}{n}$ and by \cite[(4.21.6)]{sz} the leading coefficient $L_n^{(\alpha,\beta)}$ of $P_n^{(\alpha,\beta)}$ is $L_n^{(\alpha,\beta)}=\frac{1}{2^n}\binom{2n+\alpha+\beta}{n}$. Since $|l_n|\ge c\min\{nL_{n-1}^{(\alpha+1,\beta+1)}, L_n^{(\alpha,\beta)}\}\ge c L_n^{(\alpha,\beta)}$,
$$|s_n(x)| \le c  2^n\frac{\binom{n-1+\max\{\alpha+1,\beta+1\}}{n-1}}{\binom{2n+\alpha+\beta}{n}}\le cn 2^n\frac{\Gamma(\max\{\alpha+1,\beta+1\}+n)\Gamma(\alpha+\beta+1+n)}{\Gamma(\alpha+\beta+1+2n)}$$ $$\le cn 2^n\frac{\Gamma^2(\alpha+\beta+1+n)}{\Gamma(2(\alpha+\beta)+2+2n)}\frac{\Gamma(2(\alpha+\beta)+2+2n)}{\Gamma(\alpha+\beta+1+2n)}\le c\frac{n^{\alpha+\beta+2}}{2^n}.$$
According to \eqref{ert} the estimation above ensures the result.

\medskip

\remark
Proposition \ref{c2} gives back \cite[Theorem 6.5]{bo} in 1-step Darboux transformation case, for $\alpha, \beta \ge -\frac{1}{2}$, that is recalling that $\tilde{\mu}_n=\frac{1}{n}\sum_{k=1}^{n}\delta_{x_{kn}}$ is the normalized zero-counting measure based on the
regular zeros of $\hat{P}_{n}$,
\begin{equation}\label{tmte}\tilde{\mu}_n \to \mu_e \end{equation}
in weak-star sense.

\section{Certain self-inversive polynomials} 

In this section we use the multiplication operator and the infinite matrix $M_e$ introduced in Section 3 to investigate certain self-inversive polynomials.

\noindent Let $b$ be defined by \eqref{b0}, and $\tilde{b}$ by \eqref{pperb}. Investigation of Section 3 was independent of the constant term of the polynomial $Q=\int^x \tilde{b}$. At first, as above, take $Q_0(x):= \int^x \tilde{b}$, with zero constant term. Recalling \eqref{R} we have $Q_0\hat{P}_n=\sum_{k=-L}^Lu_{n,k}\hat{P}_{n+k}$, and by \eqref{UU} $\lim_{n\to \infty}u_{n,0}=U_0$, where $U_0$ depends on $\tilde{b}$, cf. \eqref{tb}. Rearranging the equation above we have $(Q_0-U_0)\hat{P}_n=\sum_{-L\le k\le L \atop k\neq 0}u_{n,k}\hat{P}_{n+k}+ (u_{n,0}-U_0)\hat{P}_n$, thus in this section we define
$$Q(x):=\int^x \tilde{b}-U_0,$$
where $\int^x \tilde{b}$ means the primitive function without any constant term. The operator $M_e$ refers to this $Q$, that is, denoting by $M^0_e$ the operator defined by $Q_0$, $M_e=M^0_e-U_0I$. Since the asymptotics of the elements of $M_e$ can be descrided by \eqref{u3}
 $M_e$ can be decomposed to a bounded symmetric and a compact symmetric part again (cf. \eqref{ak}), that is
\begin{equation}M_{e}=M_{e,s}+M_{e,c},\end{equation}
where
$$M_{e}=\left[\begin{array}{cccccccc}u_{0,0}-U_0&u_{0,1}&\dots&\dots&u_{0,L}&0&0&\dots\\u_{1,-1}&u_{1,0}-U_0&\dots&\dots&u_{1,L-1}&u_{1,L}&0&\dots\\ \vdots&\vdots&\dots&\ddots&\vdots&\dots&\vdots&\dots\\u_{L,-L}&u_{L,-L+1}&\dots&u_{L,0}-U_0&\dots&\dots&u_{L,L}&0\\ 0&u_{L+1,-L}&\dots&\vdots&\dots&\vdots&\dots&u_{L+1,L}\\ \vdots&0&\dots&u_{L+j,-L}& \dots&\vdots&\dots&\dots\end{array}\right]$$
 and
\begin{equation}M_{e,s}=\left[\begin{array}{ccccccc}0&U_1&\dots&U_L&0&0&\dots\\U_1&0&\dots&U_{L-1}&U_L&\dots\\ \vdots&\vdots&\ddots&\vdots&\dots&\vdots&\dots\\U_{L}&\dots&0&\dots&\dots&U_{L}&\dots\\ 0&U_{L}&\vdots&\dots&\vdots&\dots&U_{L}\\ \vdots&\dots&U_{L}& \dots&\vdots&\dots&\dots\end{array}\right].\end{equation}
 Investigation of $M_{e,s}$ leads to the so-called self-inversive or palindrome polynomials, cf. \cite[Lemma 12]{mnv}. A polynomial $P$ of degree $n$ with real coefficients is self-inversive if $z^nP\left(\frac{1}{z}\right)=P(z)$. The location of zeros of self-inversive polynomials has been extensively studied, see e.g. \cite{ll}, \cite{v}, etc. Of course, the zeros of a self-inversive polynomial are symmetric with respect to the unite circle. One of the statements on location of zeros is the next one (see eg. \cite{v} ): if $P_{2m}(z)=\sum_{k=0}^{2m}a_kz^k$ is self-inversive and  $|a_m|>\sum_{0\le k \le m \atop k\neq m} |a_k|$, then $P_{2m}$ has no zeros on the unite circle. In our special case we get something similar.\\
Let us recall that $\tilde{b}(x)=\sum_{k=0}^md_kx^k$, cf. \eqref{tb}. Now define
$$P_{2L,\lambda}(z)=\sum_{k=1}^LU_k\left(z^{L+k}+z^{L-k}\right)-\lambda z^L,$$
where $U_k$ depends on $b$ see \eqref{u3}, and
$$\tilde{P}_{2L}(z)=P_{2L,\lambda}(z)+\lambda z^L.$$

\medskip

\begin{statement} $P_{2L,\lambda}(z)$ has no zeros on the unite circle if and only if
 $$\lambda \not\in \left[2\sum_{k=1}^L(-1)^kU_k,2\sum_{k=1}^LU_k\right]=\left[(-1)^L\tilde{P}_{2L}(-1),(-1)^L\tilde{P}_{2L}(1)\right].$$
\end{statement}

\medskip

 Let us consider $M_{e,s}$ as an operator on $l^2$ and on the Hardy space
 $$H^2:=\left\{f(z)=\sum_{k=0}^\infty c_kz^k : \ws \mbox{is holomorphic on} \ws |z|<1,\right.$$
 $$\left. \ws \lim_{r\to 1} f(re^{i\varphi})=f(e^{i\varphi}) \ws \mbox{a.e.} \ws \varphi \in (0,2\pi]\right\} .$$
 It is a Hilbert space under the norm $\|f\|^2=\sum_{k=0}^\infty |c_k|^2=\frac{1}{2\pi}\int_0^{2\pi}|f(e^{i\varphi})|^2d\varphi$.

 \medskip

 \begin{lemma}\label{4} $M_{e,s}-\lambda I$ has a bounded inverse if and only if $P_{2L,\lambda}(z)$ has no zeros on the unite circle .
 \end{lemma}

 \proof
 With this interpretation if $f \in H^2$,
$$(M_{e,s}f)(z)=\sum_{k=1}^LU_k(z^k+z^{-k})f(z)-\sum_{k=1}^LU_k\sum_{j=0}^{k-1}\frac{f^{(j)}(0)}{j!}z^{j-k}.$$
 Let $g\in H^2$ be arbitrary. Then considering the equation 
 $$(M_{e,s}-\lambda I)f=g$$
 we have
$$f(z)=\frac{z^Lg(z)+\sum_{k=1}^LU_k\sum_{j=0}^{k-1}\frac{f^{(j)}(0)}{j!}z^{L+j-k}}{P_{2L,\lambda}(z)}$$
\begin{equation}\label{van}=\frac{z^Lg(z)+\sum_{j=0}^{L-1}\frac{f^{(j)}(0)}{j!}\sum_{k=j}^{L-1}U_{L-k+j}z^k}{P_{2L, \lambda}(z)}.\end{equation}
 By symmetry, counting with multiplicity, $P_{2L,\lambda}$ has $l$ zeros in the unite disc and $2(L-l)$ on the unite circle. We can determine $f\in H^2$ if and only if these zeros can be compensated. That is if $\xi_m$ are the zeros of the denominator (in the closed unite disc) of multiplicity $k_m$ with $\sum_m k_m=2L-l$, then it means a system of linear equations in the numerator. According to \eqref{van} this system looks like
$$\left[\begin{array}{cccc}C_0(\xi_1)&C_1(\xi_1)&\dots&C_{L-1}(\xi_1)\\\vdots&\vdots&\dots&\vdots\\C_0(\xi_2)&C_1(\xi_2)&\dots&C_{L-1}(\xi_2)\\C'_0(\xi_2)&C'_1(\xi_2)&\dots&C'_{L-1}(\xi_2)\\ \vdots&\vdots&\dots&\vdots\\C_0^{(k_2-1)}(\xi_2)&C_1^{(k_2-1)}(\xi_2)&\dots&C_{L-1}^{(k_2-1)}(\xi_2)\\ \vdots&\vdots&\dots&\vdots\\c_{2L-l,0}&\dots&\dots&c_{2L-1,L-1}\end{array}\right]\left[\begin{array}{c}f(0)\\f'(0)\\f''(0)\\ \vdots\\ \vdots \\f^{(L-1)}(0)\end{array}\right]=\left[\begin{array}{c}h(\xi1)\\ \vdots \\h(\xi_2)\\h'(\xi_2)\\ \vdots\\ h^{(k_2-1)}(\xi_2) \\  \vdots\end{array}\right],$$
where $C_j(x)=\frac{1}{j!}\sum_{k=j}^{L-1}U_{L-k+j}x^k$, $j=0, \dots ,L-1$; $h(x)=-x^Lg(x)$, the coefficient matrix is of $2L-l \times L$. We deal with the rank of the coefficient matrix. Multiplying by an appropriate constant $k$th column and subtracting it from the $(k-1)$th one ($0<k\le L-1$) - starting the process with the last column as $C_{L-2}-(L-1)\frac{U_{L-1}}{U_L}C_{L-1}$, $C_{L-3}-(L-1)(L-2)\frac{U_{L-2}}{U_L}C_{L-1}$,... , and repeating this procedure starting with the new $(L-2)$th column, etc. it can be easily seen that finaly we get a similar matrix wich can be described on the same way as above by the modified function $\tilde{C}_j(x)=x^j$. Thus the rank of the coefficient matrix coincides with its less size. That is to get a unique solution it must be of $L\times L$ wich means that $2L-l=L$, that is $l=L$ which means that $P_{2L,\lambda}$ has $L$ zeros inside the unite disc. The unique solution of this linear system is equivalent to the existence of a (well-defined) bounded inverse of the operator. Indeed, after determination of the first $L$ coefficients of $f$, the remainder ones can be uniquely determined.  

\medskip

\begin{lemma} The spectrum of $M_{e,s}$ is $Q([-1,1])$.
 \end{lemma}

\proof
Let us recall the information on $M_{e,s}$. According to Weyl's theorem (see e.g. \cite[sec. 134]{rsz}) the essential spectrum of $M_{e,s}$ agrees with the essential spectrum of $M_{e}$. Taking into consideration that the spectrum of the multiplication operator is the closure the range of $Q$ on $[-1,1]$ (see e.g. \cite[sec. 150]{rsz}).) the essential spectrum and the spectrum of $M_{e}$ are the same. As the spectrum of $M_{e,s}$ does not contains any isolated points as well, see \cite{haw}, it also coincides with $Q([-1,1])$. \\

The fact that the essential spectrum and the spectrum of of $M_{e,s}$ are coincide in this setup can be proved as follows.\\ Let us consider $P_{2L, \lambda}(r,\varphi)$ ($z=(r\cos\varphi,r\sin\varphi)$) as a function from $\mathbb{R}^{1+2}$ to $\mathbb{R}^{2}$. If there was an isolated point of $\sigma(M_{e,s})$, then there would be a $(\lambda_0, r_0, \varphi_0)$, $r_0=1$, such that
$P_{2L, \lambda_0}(r_0,\varphi_0)=0$ and there would be a neighborhood of $\lambda_0$ such that for any $\lambda$ from this neighborhood the zeros of $P_{2L, \lambda}$ are not on the unite circle.\\
So $P_{2L, \lambda}(r,\varphi)=(p_1(\lambda,r,\varphi),p_2(\lambda,r,\varphi)$. Let us consider $\partial_2P_{2L, \lambda}=\left[\begin{array}{cc}a_{11}&a_{12}\\a_{21}&a_{22}\end{array}\right]$, where $a_{i1}=\frac{\partial}{\partial r}p_i$, $a_{i2}=\frac{\partial}{\partial \varphi}p_i$, $i=1,2$. As $\det \partial_2P_{2L, \lambda}=\frac{1}{r}(a_{22}^2+a_{12}^2)$,we can apply the implicit function theorem at $(\lambda_0,r_0,\varphi_0)$ that is there is a neighborhood of $\lambda_0$ denoted by $U$ and an arc $g$ in $U$ such that if $\lambda \in U$ then $P_{2L}(\lambda, g(\lambda))=P_{2L, \lambda}(r(\lambda),\varphi(\lambda))=0$. If $\lambda \in U$
$g'(\lambda)=\left[\begin{array}{c}r'(\lambda)\\ \varphi'(\lambda)\end{array}\right]$,\\ where $r'(\lambda)=\frac{r^{2L+1}}{a_{22}^2+a_{12}^2}\sum_{k=1}^LkU_k\cos k\varphi(r^k-r^{-k})$ (that is $r'(\lambda_0)=0$), and $\varphi'(\lambda)=\frac{-r^{2L}}{a_{22}^2+a_{12}^2}\sum_{k=1}^LkU_k\sin k\varphi(r^k+r^{-k})$. Thus the slope of the tangent line at $\lambda_0$ to $g$ is $-\cot(\varphi(\lambda_0))$, which is just the  the slope of the tangent line to the unite circle at the same point, and the curvature of $g$ at  $\lambda_0$ is $1$. Considering the symmetry of the zeros with respect to the unite circle $g$ has to coincide locally with the unite circle, which means that $\lambda_0$ cannot be isolated.

\medskip

\proof (of the Statement)
Notice, that $b\ge 0$ on $[-1,1]$, thus $Q$ is increasing here and $Q(-1)<Q(1)$. That is to prove the statement it is enough to compute these two values. In view of \eqref{u3}
$$2\sum_{k=1}^L(\pm 1)^kU_k$$ $$=2\sum_{l=1}^{\left[\frac{L}{2}\right]}\sum_{p=l}^{\left[\frac{L}{2}\right]}\frac{d_{2p-1}}{2p}\binom{2p}{p-l}\frac{1}{2^{2p}}\pm 2\sum_{l=0}^{\left[\frac{L-1}{2}\right]}\sum_{p=l}^{\left[\frac{L-1}{2}\right]}\frac{d_{2p}}{2p+1}\binom{2p+1}{p-l}\frac{1}{2^{2p+1}}=:S.$$
Changing the order of summation and by the definition of $Q$
$$S=\sum_{p=l}^{\left[\frac{L}{2}\right]}\frac{d_{2p-1}}{2p}\left(1-\frac{1}{2^{2p}}\binom{2p}{p}\right)\pm \sum_{p=l}^{\left[\frac{L-1}{2}\right]}\frac{d_{2p}}{2p+1}=\sum_{k=1}^L(\pm 1)^k\frac{d_{k-1}}{k}-U_0=Q(\pm 1).$$
Thus the spectrum of the operator $[M]_{e,s}$ is $\left[\sum_{k=1}^L(-1)^kU_k,2\sum_{k=1}^LU_k\right]$, which together with Lemma \ref{4} proves the statement.

\medskip

\noindent \small{Department of Analysis, \newline
Budapest University of Technology and Economics}\newline
\small{ g.horvath.agota@renyi.hu}

\end{document}